\documentclass{eccomas2016}

\usepackage{amsmath}
\usepackage{amssymb}
\usepackage{subcaption}

\title{HYBRID RIEMANN SOLVERS FOR LARGE SYSTEMS OF CONSERVATION LAWS}

\author{Birte Schmidtmann$^1$, Mariia Astrakhantceva$^2$, and Manuel Torrilhon$^1$}

\heading{Birte Schmidtmann, Mariia Astrakhantceva, and Manuel Torrilhon}

\address{$^1$ MathCCES, RWTH Aachen University \\
  Schinkelstr. 2, 52062 Aachen, Germany\\
  e-mail: \{schmidtmann, mt\}@mathcces.rwth-aachen.de \and
  $^2$
  Lucerne School of Business\\
  Zentralstrasse 9, Postfach 2940, 6002 Luzern, Switzerland\\
  e-mail: mariia.astrakhantceva@hslu.ch}


\abstract{In this paper we present a new family of approximate Riemann solvers for the numerical approximation of solutions of hyperbolic conservation laws. They are approximate, also referred to as incomplete, in the sense that the solvers avoid computing the characteristic decomposition of the flux Jacobian. Instead, they require only an estimate of the globally fastest wave speeds in both directions. Thus, this family of solvers is particularly efficient for large systems of conservation laws, i.e. with many different propagation speeds, and when no explicit expression for the eigensystem is available. 
Even though only fastest wave speeds are needed as input values, the new family of Riemann solvers reproduces all waves with less dissipation than HLL, which has the same prerequisites, requiring only one additional flux evaluation.}

\begin{document}
\section{INTRODUCTION}
In the finite volume method, integrating conservation laws over a control volume leads to a formulation which requires the evaluation of the flux function at cell interfaces. Because the exact information is not available, a sequence of local Riemann problems needs to be solved \cite{LeVeque1992}. The initial states for these problems are typically given by the left and right adjacent cell values, $(u_i^n, u_{i+1}^n)$ for each $i$ \cite{ToroRiemannSolvers}. These local Riemann problems have to be solved many times for finding the numerical solution, therefore, the Riemann solver is a building block of the finite volume method. Note that the numerical flux function at interfaces also appears in the Discontinuous Galerkin (DG) formulation. Therefore, Riemann solvers are also needed in DG.

Over the last decades, many different Riemann solvers have been developed, see e.g. \cite{ToroRiemannSolvers}. The challenge is that the solver needs to be computationally efficient and easy to implement. At the same time, it needs to yield accurate results which do not create artificial oscillations. The latter is ensured by requiring the solver to be monotone.

Riemann solvers can be classified into complete and incomplete schemes, depending on whether all present characteristic fields are considered in the model or not. According to this classification, the upwind scheme and Roe's scheme \cite{Roe1981}, respectively, are complete schemes. They yield good, monotone results, however, an evaluation of the eigensystem of the flux Jacobian is needed. This characteristic decomposition might be expensive to compute, especially for large systems, and in some cases, an analytic expression is not available at all. However, if possible, using Roe's scheme typically yields the best resolution of the Riemann wave fan, since all waves are well-resolved. Therefore, this scheme can be considered as the optimum, also taking into account that less dissipation than upwinding leads to non-monotone solutions. In order to keep its high resolution and at the same time reducing the computational cost, there have been many attempts to approximate the upwind scheme without solving the eigenvalue problem, see e.g. \cite{Torrilhon2012, CastroGallardoMarquina2014} and references therein.	
	
	In this article, we are interested in approximate Riemann solvers. Their advantage is the easy implementation and the requirement of few characteristic information. However, incomplete Riemann solvers contain more dissipation than the upwind scheme and thus, yield lower resolution, especially of slow waves. Nevertheless, in many test cases, these Riemann solvers may be sufficient to obtain good results, especially if the system contains only fast waves. 
\section{FINITE VOLUME METHOD}\label{sec:FV}

We consider a hyperbolic system of conservation laws in one space dimension of the form
\begin{align}
	\label{eq:consLaw}
	\partial_t\,U(x,t)+\partial_x f(U(x,t))&=0\quad\text{in}\;\mathbb{R}\times\mathbb{R}^+ \\
	U(x,0)&=U_0(x)
\end{align}
where the unknown variable vector is given by $U:\mathbb{R}\times\mathbb{R}^+\rightarrow\mathbb{R}^N$, with some initial conditions $U_0(x)$ and a flux function $f:\mathbb{R}^N\rightarrow\mathbb{R}^N$, such that the Jacobian $A(U)=D f(U)$ only has real eigenvalues. In finite volume methods, the domain is subdivided into cells $C_i=[x_{i-1/2}, x_{i+1/2}]$ and Equation \eqref{eq:consLaw} is integrated over a cell $C_i$. The cell average at time $t^n$, given by
\begin{align}
	\frac{1}{\Delta x}\int_{C_i} U(x, t^n) dx,
\end{align}
is approximated by $\bar U_i^n$. The update of the approximation at time $t^{n+1}=t^n+\Delta t$ reads \cite{LeVeque2002}
\begin{align}
	\bar U_i^{n+1} = \bar U_i^{n} - \frac{\Delta t}{\Delta x} \left( \hat f_{i+1/2}-\hat f_{i-1/2} \right)
\end{align}
and contains the numerical flux function $\hat f_{i+1/2}=\hat f(U_{i+1/2}^{-}, U_{i+1/2}^{+})$. This function takes as input values the left and right limiting values of the variable vector $U$ at the cell interface $i+1/2$. 
 
 The choice of the numerical flux function - or Riemann solver - of the method, determines the mathematical properties of the scheme, such as accuracy and monotonicity. Thus, the choice of the solver is crucial for the resulting scheme. 
 
 In this work, we consider the general formulation of the numerical flux, given by
\begin{align}
	\label{eq:numFlux}
	\hat f (U_L, U_R) = \frac{1}{2}\left( f(U_L) + f(U_R)\right) + \frac{1}{2}D(U_L, U_R) \left(U_L - U_R\right)
\end{align}
with a dissipation matrix $D$ which depends on the left and right states $U_L$ and $U_R$, respectively. $D$ also depends on characteristics of the flux Jacobian $A(U)$. 

Let us denote by $\lambda_\text{min}(U)$ and $\lambda_\text{max}(U)$ the minimal and maximal eigenvalue of $A(U)$. The spectral radius at $U$, i.e. the maximum absolute characteristic speed, is given by \\ $\bar \lambda = \max\{|\lambda_\text{min}|, |\lambda_\text{max}|\}$.

\section{REVIEW OF RIEMANN SOLVERS}\label{sec:riemannsolvers}
In this section we recall some well-known Riemann solvers which are necessary for the development of the new family of hybrid Riemann Solvers. First, we note that the dissipation matrix $D$ completely dictates the numerical flux function \eqref{eq:numFlux} and hence the numerical scheme. Therefore, we break down the following discussion to comparing the dissipation matrices of the schemes. Furthermore, for the sake of comparability, we introduce the notation of the dimensionless scalar dissipation function $d(\nu)$:

The flux Jacobian $A(U)$ can be diagonalized as $A(U)=T(U)\,\Lambda(U)\,T(U)^{-1}$, where $\Lambda(U)$ is the eigenvalue matrix $\Lambda(U)=\text{diag}(\lambda_1(U),\ldots,\lambda_N(U))$, with $\lambda_1<\lambda_2<\ldots<\lambda_N$, and $T(U)$ the corresponding eigenvector matrix. Since the dissipation matrix $D$ is a function of the flux Jacobian $A$, it can be shown that
	\begin{align}
		\frac{\Delta t}{\Delta x}D(A)=T^{-1}\,\frac{\Delta t}{\Delta x}\,\text{diag}(\tilde d(\lambda_1), \ldots, \tilde d(\lambda_N))\,T=T^{-1}\,\text{diag}(d(\nu_1), \ldots, d(\nu_N))\,T,
		\label{eq:trafoDd}
	\end{align}
	holds. Here, $\nu_i=\lambda_i\Delta t/\Delta x$ and $d(\nu)$ is the dimensionless scalar dissipation function. The eigenvector matrix $T$ is the same, independent of $D$.
	
	We will compare the dissipation functions of different schemes in a $\nu-d(\nu)$-plot at the end of this section.
\newline\newline
For a linear system, where we can write $f(U)=A\,U$ with some constant matrix $A \in \mathbb{R}^{N \times N}$, the dissipation matrix of the upwind scheme reads
	\begin{align}
		D_\text{up}=|A|\quad\leftrightarrow\quad d_\text{up}(\nu)=|\nu|.
	\end{align}
	For non-linear systems with general flux function, the scheme has been extended by Roe \cite{Roe1981}. The dissipation matrix of Roe's scheme is given by
	\begin{align}
		D_\text{Roe}=|\tilde A|\quad\leftrightarrow\quad d_\text{Roe}(\nu)=|\nu|,
	\end{align}
	with a so-called Roe Matrix $\tilde A$, satisfying $\tilde A(U_L, U_R) (U_L-U_R) = F(U_L) - F(U_R)$. The upwind Godunov solver and its non-linear extension, the Roe solver are complete Riemann solvers. Now follows a list of incomplete solvers with decreasing dissipation.
	
	The dissipation function of the monotone Lax-Friedrichs scheme is
	\begin{align}
		D_\text{LF}=\frac{\Delta x}{\Delta t}\,I\quad\leftrightarrow\quad d_\text{LF}(\nu)=1.
	\end{align}	
	A solver which decreases the dissipation is the Rusanov scheme, also referred to as local Lax-Friedrichs scheme. It takes into account the globally fastest eigenvalue of the system:
	\begin{align}
		D_\text{LLF}=\max\{\bar{\lambda}(U_L), \bar{\lambda}(U_R)\}\,I\quad\leftrightarrow\quad d_\text{LLF}(\nu)=\max(|\nu_{\min}|,|\nu_{\max}|).
	\end{align}	
	Harten, Lax and van Leer \cite{HLL1983} further decreased the dissipation, especially of slow waves, by considering the fastest and slowest waves of the system:
\begin{subequations}
	\begin{align}
		D_\text{HLL} &= \frac{|\lambda_L|-|\lambda_R|}{\lambda_L-\lambda_R}\,\tilde{A} - \frac{|\lambda_L|\,\lambda_R - |\lambda_R|\,\lambda_L}{\lambda_L-\lambda_R}\,I\\
		\leftrightarrow d_\text{HLL} &= \frac{|\nu_L|-|\nu_R|}{\nu_L-\nu_R}\,\nu - \frac{|\nu_L|\,\nu_R - |\nu_R|\,\nu_L}{\nu_L-\nu_R}.
	\end{align}
\end{subequations}
Here, $\lambda_L=\lambda_\text{min}(U_L)$, $\lambda_R=\lambda_\text{max}(U_R)$, and $\nu_{L/R}=\lambda_{L/R}\Delta t/\Delta x$.

A scheme which further reduces the dissipation is the Lax-Wendroff scheme,
	\begin{align}
		D_\text{LW}=\frac{\Delta t}{\Delta x}\,\tilde A^2 \quad\leftrightarrow\quad d_\text{LW}(\nu)=\nu^2.
	\end{align}	
	Discontinuities are approximated with steep gradients using the Lax-Wendroff scheme, however, the method is known to cause oscillations at discontinuities because it is non-monotone \cite{LeVeque1992}, see also Fig.~\ref{fig:differentSolvers}. 
\section{A FAMILY OF NEW HYBRID RIEMANN SOLVERS}\label{sec:HLLXomega}
%
We want to construct Riemann solvers which require as few information as HLL but are less dissipative. This is advantageous for the resolution, especially for slow waves, i.e. $\lambda \approx 0$. Demanding the same input values as HLL, we require the knowledge or an estimate of the globally slowest and fastest characteristic waves of the system. Only one additional flux evaluation shall be required.

%
\subsection{$P_2$ - conditions}\label{subsec:HLLX}
\begin{figure}[t]
	\centering %
	\includegraphics[width=0.5\textwidth]{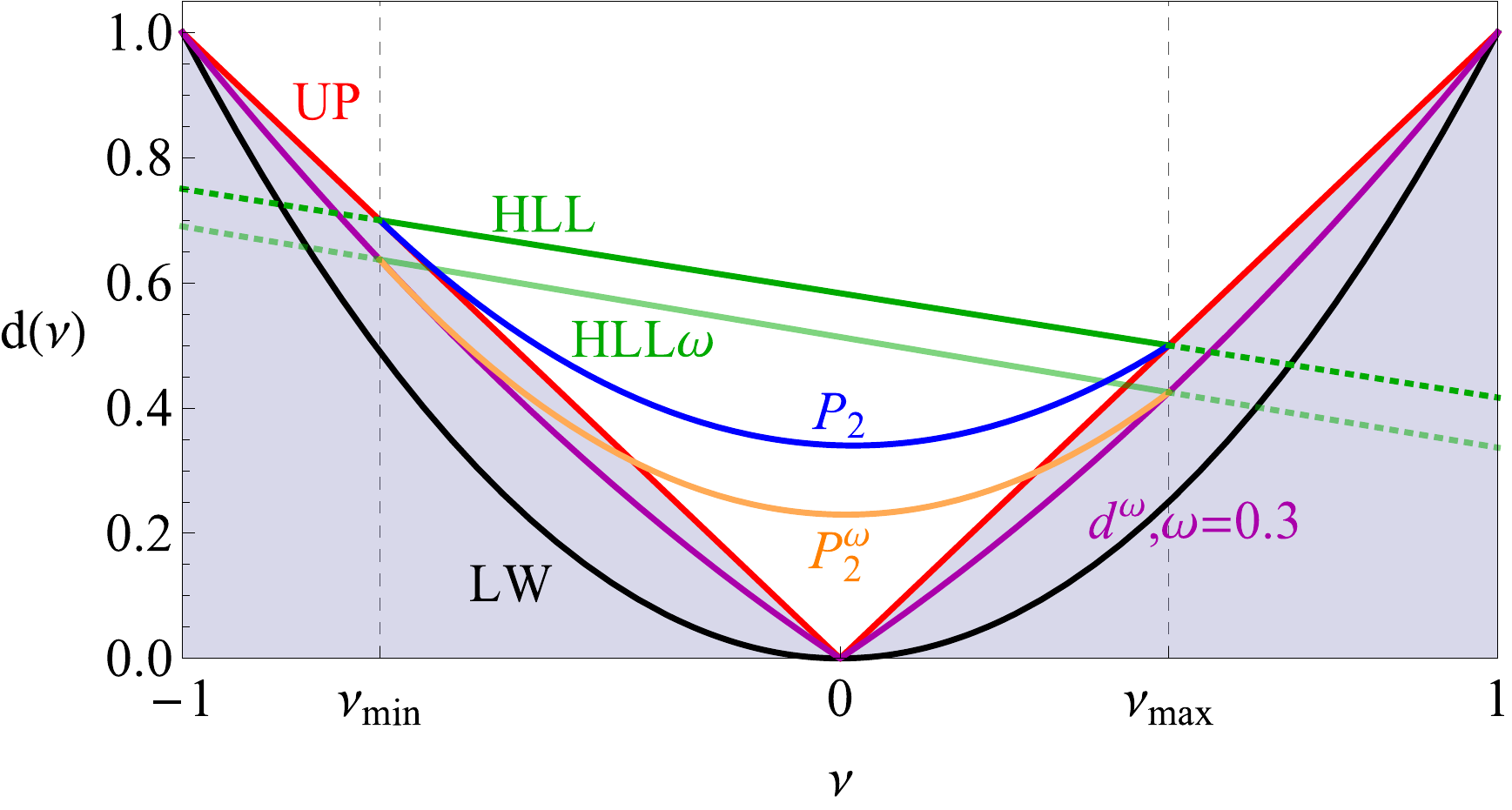}
	\caption{Scalar, non-dimensional dissipation functions of different Riemann solvers.}
	\label{fig:differentSolvers}
\end{figure}

%
We now pick up three requirements which have been proposed by Degond et al. \cite{DegondPeyrardRussoVilledieu1999}. The resulting monotone Riemann solver, named $P_2$, is based on a quadratic dissipation function, fully determined by 
\begin{subequations}
	\begin{align}
		d_{P_2}(\nu_\text{min})&=d_\text{up}(\nu_\text{min})=|\nu_\text{min}|,\\
		d_{P_2}(\nu_\text{max})&=d_\text{up}(\nu_\text{max})=|\nu_\text{max}|, \,\text{and}\\
		d_{P_2}'(\bar\nu)&=d_\text{up}'(\bar\nu)=\text{sign}(\bar\nu), \quad\bar\nu=\begin{cases}\nu_\text{max},\quad\text{if}\;\;|\nu_\text{max}|\geq|\nu_\text{min}|\\
			\nu_\text{min},\quad\text{if}\;\;|\nu_\text{min}|>|\nu_\text{max}|.	\end{cases}
	\end{align}
\end{subequations}	
This function automatically fulfills $d_{P_2}(\nu)\geq |\nu|$ for $\nu \in [\nu_\text{min}, \nu_\text{max}]$ which means that it is monotone in this region. This can also be seen in Fig.~\ref{fig:differentSolvers}, which additionally shows that $P_2$ is less dissipative than HLL, especially for $\nu\approx 0$. The dissipation function $d_{P_2}(\nu)$ can be written in the simple form
\begin{align}
	d_{P_2}(\nu) =d_\text{HLL}(\nu) + \alpha  (\nu-\nu_\text{min})(\nu-\nu_\text{max}),
	\label{eq:dHLLX}
\end{align}
with
\begin{align}
	\alpha = \frac{\nu_\text{max}- \nu_\text{min} - \big| |\nu_\text{max}| - |\nu_\text{min}| \big| }{(\nu_\text{max}-\nu_\text{min})^2}.
	\label{eq:alpha}
\end{align}

%
\subsection{$P_2^\omega$ - Beyond Monotonicity}\label{subsec:HLLXomega}
%
The main idea of $P_2^\omega$ is to construct a quadratic dissipation function, in the form of $P_2$. This new function shall be closer to the absolute value function, i.e. the upwind scheme, for all waves $\lambda_i$ (and thus for all $\nu_i$) of the hyperbolic system. Since we do not want to increase neither the number of input information, nor the number of flux evaluations, we lower the dissipation function by a certain amount. This amount is described by a parameter $\omega\in [0,1]$, which determines the monotonicity behavior of the solver. For $\omega=0$ we recover the monotone $P_2$ solver, and for $\omega=1$, the non-monotone Lax-Wendroff solver. All intermediate members of the $P_2^\omega$ family are slightly non-monotone for a certain range of waves. However, we show in this section that under some mild assumptions, the results do not show spurious oscillations.
%
\subsubsection{Monotonicity Study}
%
Before we introduce this family of Riemann solvers, we state and validate some observations made by Torrilhon \cite{Torrilhon2012}. Firstly, it was perceived that the MUSTA fluxes introduced by Toro \cite{Toro2006} slightly drop below the upwind flux, which means that they do not fully lie in the monotonicity preserving region. Thus, as expected, the numerical solutions obtained with MUSTA fluxes show some non-monotone behavior. However, this behavior is far from the oscillations created by the Lax-Wendroff scheme. Additionally, the observed oscillations of MUSTA solutions decay in time and disappear after a certain number of time steps, cf. \cite[Fig. 5, p. A2084]{Torrilhon2012}. These results are essentially independent of the grid size.

These interesting results were observed for a dissipation function which slightly drops below the absolute value function. Let us introduce the dissipation function $d_\omega(\nu)$, which is the weighted average of the dissipation functions of the monotone upwind scheme and the non-monotone Lax-Wendroff scheme,
\begin{align}	
	\label{eq:dOmega}
	d_\omega(\nu, \omega) &=\omega\ d_{\text{LW}}(\nu)+(1-\omega)\ d_{\text{up}}(\nu)\;\; \omega\in [0,1].
\end{align}	
For $\omega=0$ we recover the monotone upwind scheme $d_{\omega=0}(\nu) = d_{UP}(\nu)$ and for $\omega=1, d_{\omega=1}(\nu) = d_{LW}(\nu)$ holds true. Fig.~\ref{fig:differentSolvers} shows $d_\omega(\nu, \omega)$ for $\omega=0.3$. 

The aim of this section is to study the monotonicity behavior of $d_\omega(\nu)$ and produce similar effects as studied in \cite{Torrilhon2012}. Let us therefore investigate the solutions of the numerical flux function Eq.~\eqref{eq:numFlux} with Eq.~\eqref{eq:dOmega} for different values of $\omega$. The test problem is the scalar transport equation $u_t + u_x=0$ with initial condition $u_0(x)=\text{sign}(x)$ on the interval $[-1, 1]$. The jump evolves with time on a grid with $n=200$ grid cells until $T_\text{end}=0.25$. The Courant number is set to $CFL=0.5$, which shows the maximal deviation of Lax-Wendroff from Upwind. 	
	\begin{figure}[t]
		\centering
		\begin{subfigure}{.4\textwidth}
			\includegraphics[width=\textwidth]{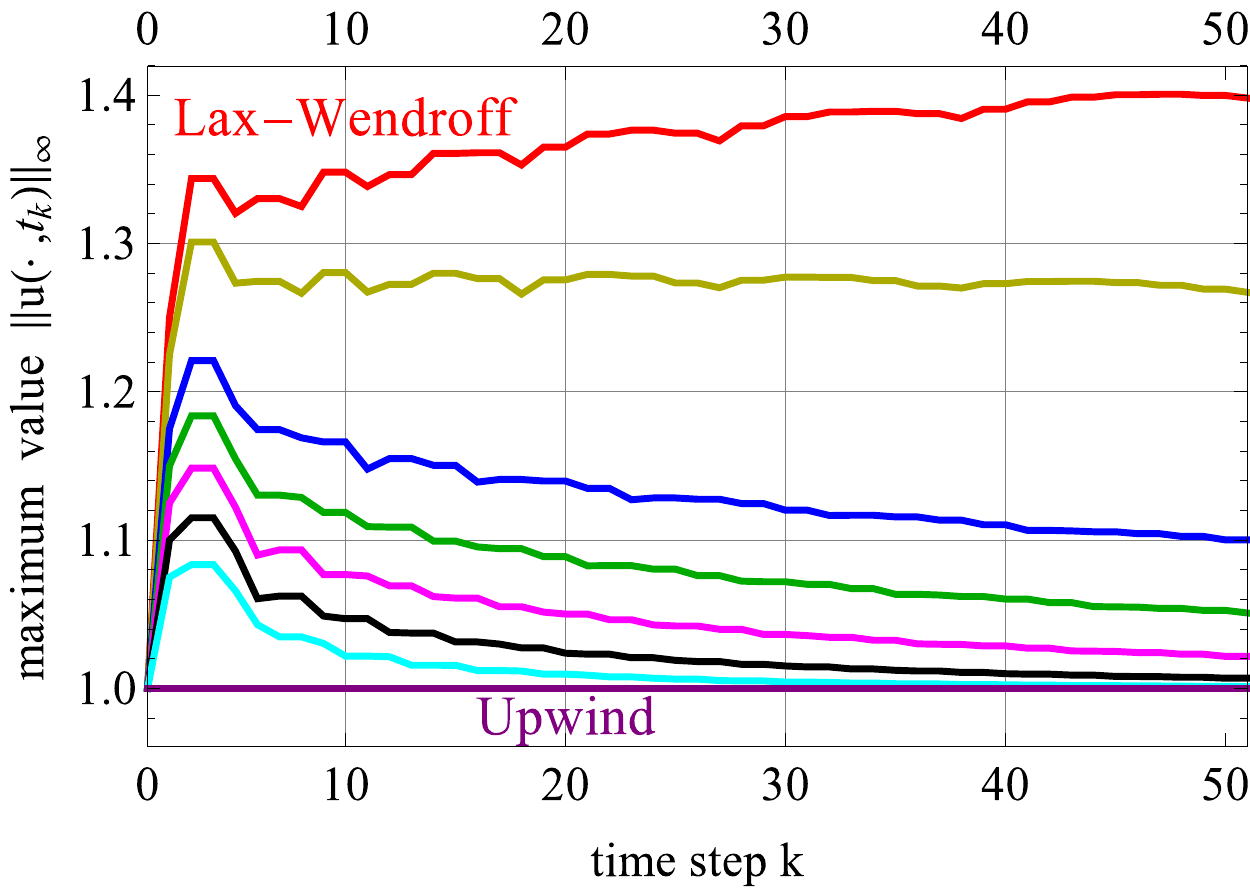}	
			\caption{Maximum value of the solution $u$ as a function of the number of time steps.}
			\label{fig:maxDomega}
		\end{subfigure}	
		\hfill
		\begin{subfigure}{.56\textwidth}
			\includegraphics[width=\textwidth]{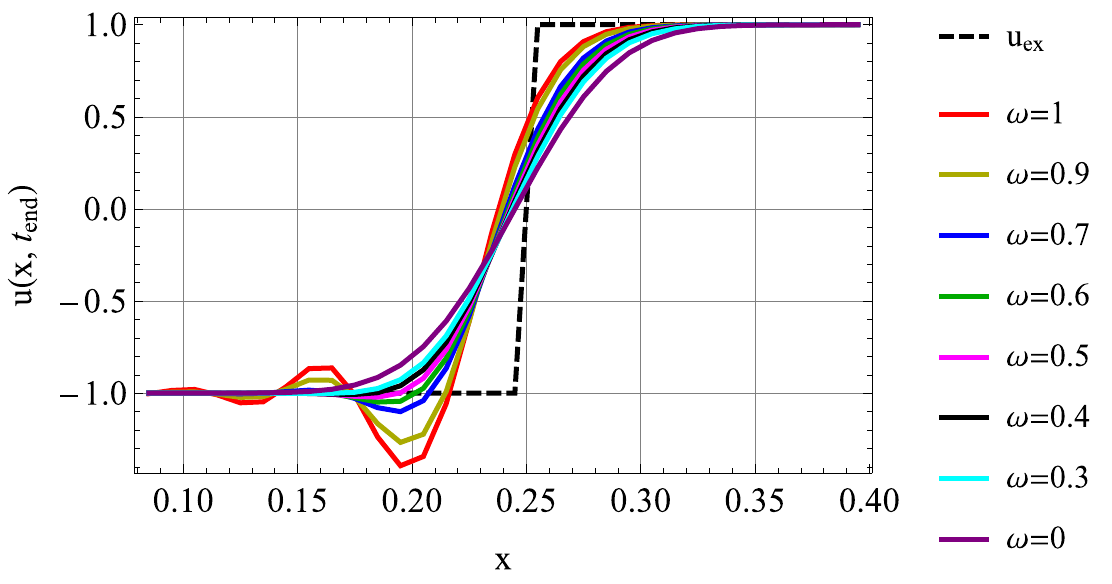}
			\caption{Zoom of the solution for different $\omega$.}
			\label{fig:uOmegas}
		\end{subfigure}	
		\caption{Test case with initial condition sign$(x)$ on $x\in [-1, 1]$ with $n=200$ grid cells, $\,CFL= 0.5,$ and end-time $T_{\text{end}}=0.25$, which corresponds to $50$ time steps.}
	\end{figure}

	The numerical results for all tested values of $\omega$ are shown in Fig.~\ref{fig:uOmegas}. It can be easily seen, that $\omega=0$ and $\omega=1$ correspond to the upwind and the Lax-Wendroff schemes. That is, for $\omega=1$ we can observe the well-known oscillations. As $\omega$ decreases, the oscillations also decrease and for $\omega \leq 0.4$ no oscillations can be seen anymore. This can also be seen in Fig.~\ref{fig:maxDomega}, which shows the maximum value of the solutions $u$ depending on the number of time steps. Here, 50 time steps correspond to $T_{\text{end}}=0.25$. This figure shows that the oscillations, which appear in the Lax-Wendroff scheme do not decrease over time. However, as soon as the monotone upwind scheme has enough weight, the oscillations start decreasing over time. This phenomenon can be explained by looking at the modified equations \cite{LeVeque1992} of upwind, Lax-Wendroff and the their weighted average. In the modified equation of the latter, the diffusive upwind term remains the dominant term. This avoids the creation of oscillations caused by the Lax-Wendroff dispersion term, when this term has enough weight. More details on the modified equations will be presented in future work.

	$d_\omega(\nu, \omega)$ could also be used as a limiter function, increasing $\omega$ in smooth parts of the solution and setting it to $0$ (i.e. recovering the upwind solver) when discontinuities are encountered.
%
\subsubsection{HLL$\omega$}
%
Based on the findings above, let us define a Riemann solver which is a modification of HLL with less dissipation. We shall call this solver HLL$\omega$. In the same way HLL is constructed \cite{HLL1983, ToroRiemannSolvers}, the dissipation function of this solver is of the form $d_{\text{HLL}\omega}(\nu) = b_0 + b_1 \nu$. Here, the coefficients depend on $\omega$, i.e. $b_0=b_0(\omega), b_1=b_1(\omega)$. Additionally, instead of intersecting with the absolute value function at $\nu_\text{min}$ and $\nu_\text{max}$, HLL$\omega$ fulfills the following constraints:
	\begin{align}
		d_{\text{HLL}\omega}(\nu_{\min})=d_\omega(\nu_{\min}), \quad d_{\text{HLL}\omega}(\nu_{\max})=d_\omega(\nu_{\max}).
	\end{align}
	The dissipation function is shown in Fig.~\ref{fig:differentSolvers}, where it is well-visible, that HLL$\omega$ is less dissipative than HLL and is non-monotone for some wave speeds.
%
\subsubsection{$P_2^\omega$}
%
Now we can come back to the aim of this section, the construction of a new family of approximate Riemann solvers - $P_2^\omega$. The dissipation functions of these solvers are similar to $d_{P_2}$, only closer to the absolute value function for all emerging wave speeds of the system. The dissipation functions are given by
\begin{align}
	d_{P_2^\omega}(\nu) =d_{\text{HLL}\omega}(\nu) + \beta (\nu-\nu_\text{min})(\nu-\nu_\text{max}),
\end{align}
with 
\begin{align}
	\beta = \omega + (1-\omega)\frac{\nu_\text{max}-\nu_\text{min}-\big| |\nu_\text{max}| - |\nu_\text{min}| \big|}{(\nu_\text{max}-\nu_\text{min})^2}.
	\label{eq:beta}
\end{align}
Note that $\beta = \omega + (1-\omega)\alpha$, with the $P_2$ coefficient $\alpha$ \eqref{eq:alpha}. Thus, it is easy to verify that for $\omega=0$, the monotone $P_2$ solver is recovered. 

It can be seen in Fig.~\ref{fig:differentSolvers} that $P_2^\omega$ is less dissipative than HLL and $P_2$. However, it does not fully lie in the monotonicity preserving region, thus, one would expect some non-monotone behavior. Nevertheless, we observe that oscillations appearing close to discontinuities disappear after a certain number of time steps. Thus, the final result obtained with $P_2^\omega$ is non-oscillatory.
	
The choice of $\omega$ remains problem-dependent. However, $\omega\leq 0.5$ turned out to be a good choice and will be used in the following.
%
\section{NUMERICAL RESULTS - APPLICATION TO IDEAL MHD}\label{sec:numericalresults}
%
We will apply the new family of Riemann solvers to the equations of ideal magnetohydrodynamics (MHD) and compare it to other solvers. Ideal MHD describes the flow of plasma, assuming infinite electrical resistivity. The equations of ideal MHD in one-dimensional processes read
\begin{align}
\label{eq:idealMHD}
	\partial_t \begin{pmatrix}
	\rho\\ \rho v_x\\ \rho\bold{v_t}\\ \bold{B_t}\\ E
\end{pmatrix}	 
+ \partial_x \begin{pmatrix}
	\rho v_x\\ \rho v_x^2+p+\tfrac{1}{2}\bold{B_t^2}\\ \rho v_x\bold{v_t}-B_x\bold{B_t}\\ v_x\bold{B_t}-B_x\bold{v_t}
	\\ (E+p+\tfrac{1}{2}\bold{B_t^2})v_x - B_x\bold{B_t}\cdot\bold{v_t}
\end{pmatrix}
 = 0
\end{align}
with density $\rho$, normal and tangential velocities $v_x$, and $\bold{v_t}=(v_y, v_z)$, respectively. The normal magnetic field $B_x$ is constant in the one-dimensional case, the tangential magnetic field is $\bold{B_t}=(B_y, B_z)$ and the energy $E$ is given by
\begin{align}
	\label{eq:energy}
	E=\frac{1}{\gamma-1}p + \frac{1}{2}\rho (v_x^2 + \bold{v_t^2})+\frac{1}{2}\bold{B_t^2}
\end{align}
in terms of the pressure $p$. The adiabatic constant $\gamma$ is set to $5/3$. Since the normal magnetic field $B_x$ is constant, system \eqref{eq:idealMHD} contains seven equations for the seven unknowns, $\bold{U}=(\rho, v_x,\bold{v_t}, p, \bold{B_t})$, exhibiting seven characteristic velocities, and therefore can be considered as a large system of conservation laws.

Let us consider the Riemann problem given by
\begin{align}
	\label{eq:ICmhd}
	\bold{U}_0(x)= \begin{cases} \bold{U}_L=(3,0,(0,0),3,(1,1))\;\;\text{if}\;x<0\\
						 \bold{U}_R=(1,0,(0,0),1,(\cos(1.5), \sin(1.5)))\;\;\text{if}\;x\geq 0,
			\end{cases}			 
\end{align}
\begin{figure}[t]
	\centering
	\begin{subfigure}{.469\textwidth}
		\includegraphics[width=\textwidth]{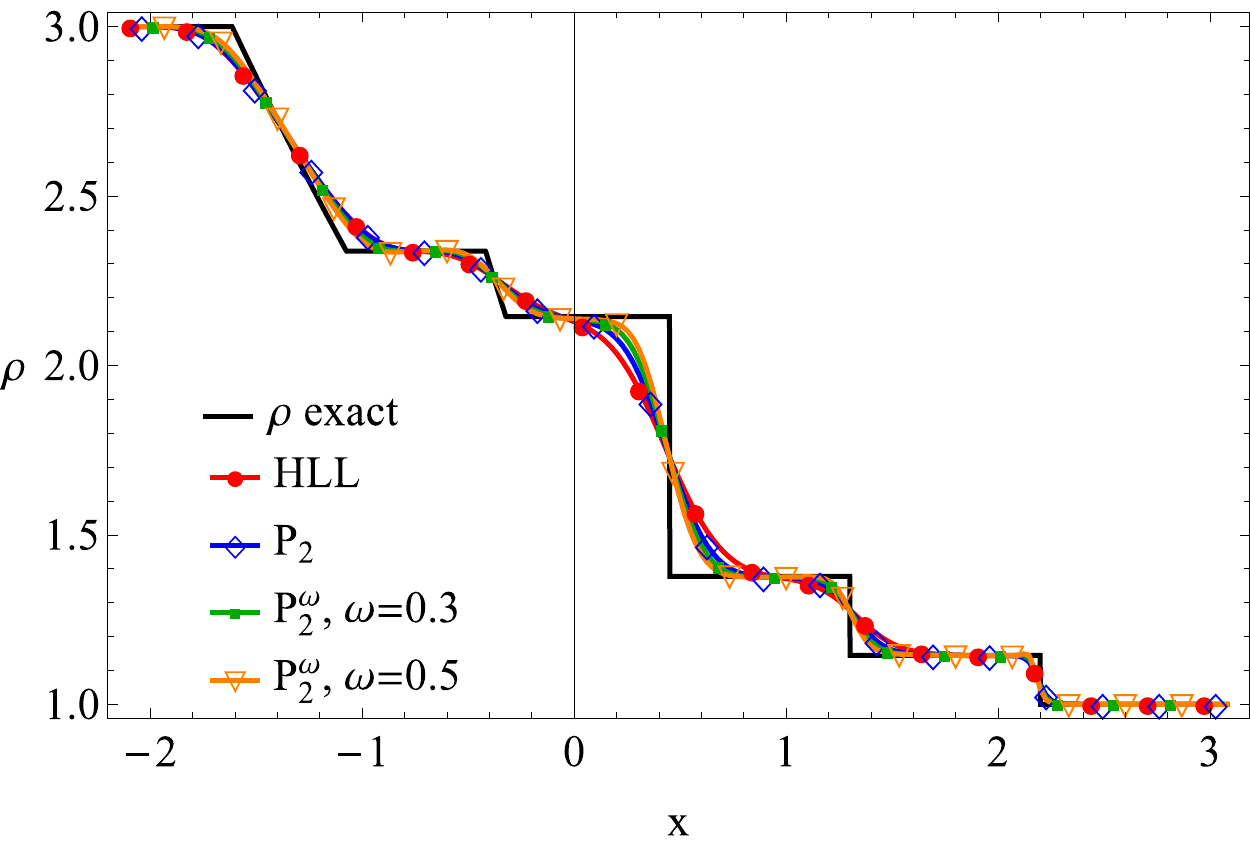}
		\caption{Density profile for different solvers.}
		\label{fig:density}
	\end{subfigure}	
	\hfill
	\begin{subfigure}{.47\textwidth}
		\includegraphics[width=\textwidth]{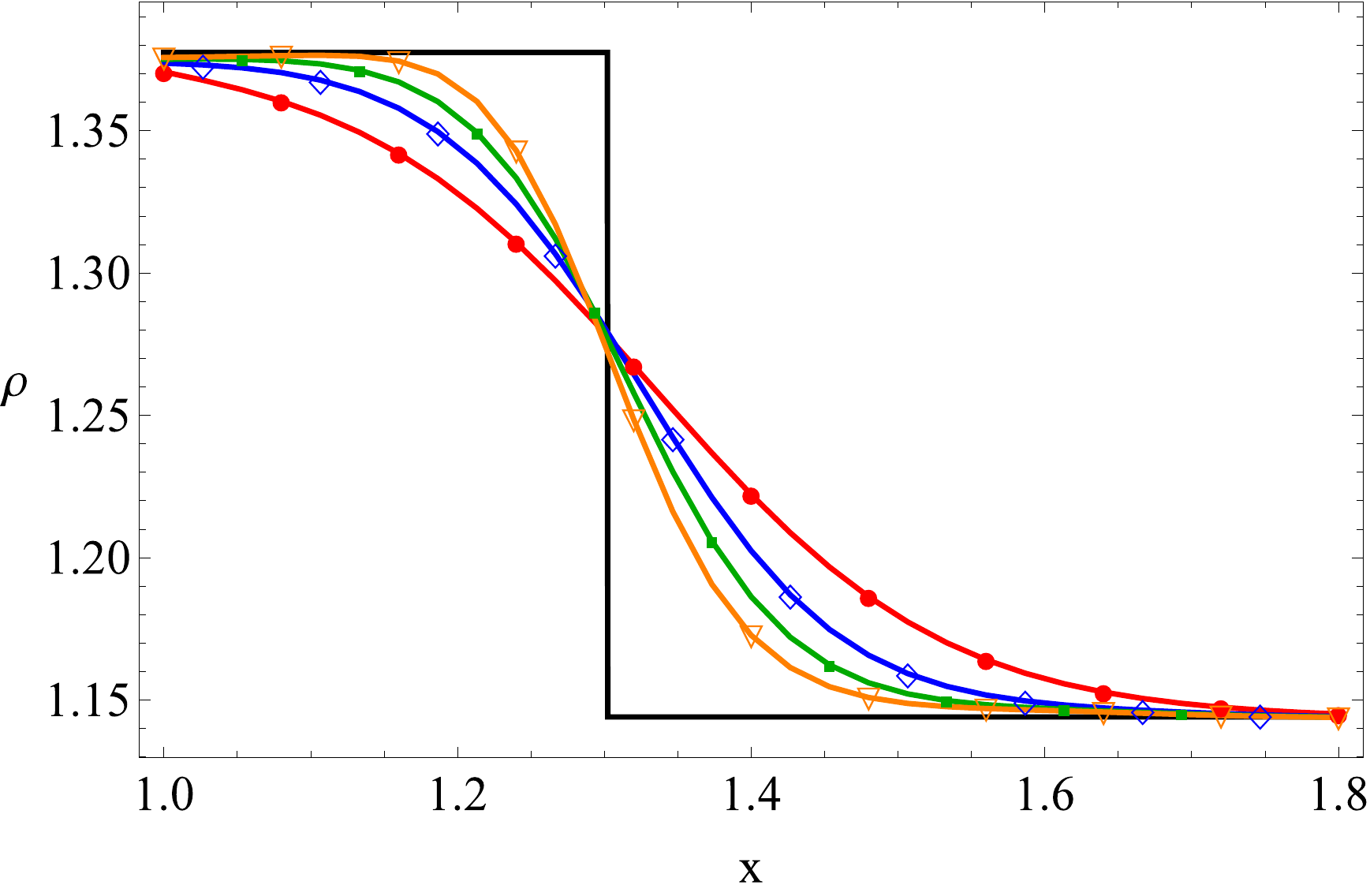}
	\caption{Zoom of the slow shock around $x=1.3$.}
	\label{fig:densityZoom}	
	\end{subfigure}	
	\caption{Solution of the ideal MHD equations \eqref{eq:idealMHD} with initial conditions \eqref{eq:ICmhd} in the domain $x\in [-4, 4]$ with $n=300$ grid cells, $\,CFL= 0.9,\, \Delta t=0.01$, and end time $T_{\text{end}}=1.0$.}
	\label{fig:solMHD}
\end{figure}
and $B_x=1.5$. The computational domain is $[-4, 4]$, and the solutions depicted in Fig.~\ref{fig:solMHD} have been obtained with $N=300$ grid cells, CFL number $0.9$, $\Delta t=0.01$, and end time $T_\text{end}=1.0$. Fig.~\ref{fig:density} shows the solution of the density profile obtained with HLL, $P_2$ and $P_2^\omega$, $\omega= 0.3, 0.5$. The exact solution has been obtained by \cite{exactMHDsolver}. Fig.~\ref{fig:densityZoom} shows a zoom of the slow right moving shock at $x=1.3$, where the resolution of all solvers can be nicely compared. It can be stated that $P_2$ increases the resolution compared to HLL, causing a steeper gradient. $P_2^\omega$ further increases the steepness of the gradient for increasing $\omega$, due to decreasing dissipation. This effect is present in an even stronger form at the contact discontinuity, which corresponds to a slower wave. At the fast shock, the differences between the solutions of the four solvers are less significant, because this discontinuity relates to a larger wave speed $\lambda$. This observation corresponds well to Fig.~\ref{fig:differentSolvers}, where we showed that the differences of the discussed dissipation functions are larger for slower waves.
%

\section{CONCLUSIONS}\label{sec:conclusions}

This paper presented a family of approximate hybrid Riemann solvers, $P_2^\omega$, for non-linear hyperbolic systems of conservation laws.
The solvers do not require the characteristic decomposition of the flux Jacobian, only an estimate of the maximal propagation speeds in both directions is needed. The family of solvers contains a parameter $\omega$ which orders the solvers from fully-monotone to fully non-monotone. The intermediate solvers contain monotone as well as non-monotone parts. We showed that these intermediate family members, even though containing non-monotone parts for certain wave speeds, do not lead to oscillatory solutions. 

Extremely slow waves and stationary waves will still be approximated with higher dissipation than the upwind scheme, however, the computational cost of the new solvers is lower. Compared to solvers with similar prerequisites, the new Riemann solvers are able to rigorously decrease the dissipation of the scheme.

The numerical examples underline the excellent performance of the new family of solvers with respect to other solvers.

\end{document}